\documentclass[12pt]{amsart}
\calclayout

\let\le\undefined
\DeclareMathSymbol{\le}{\mathrel}{AMSa}{"36}      
\let\ge\undefined
\DeclareMathSymbol{\ge}{\mathrel}{AMSa}{"3E}      
\let\k\undefined
\DeclareMathSymbol{\k}{\mathalpha}{AMSb}{"7C}
\DeclareMathSymbol{\lefttimes}{\mathbin}{AMSa}{"68}
\DeclareMathSymbol{\righttimes}{\mathbin}{AMSa}{"69}
\let\empty\undefined
\DeclareMathSymbol{\empty}{\mathord}{AMSb}{"3F}

\newcommand{\ds}{\dots}
\newcommand{\dsm}{\dotsm}

\newcommand{\ot}{\otimes}
\newcommand{\rk}{\operatorname{rk}}
\newcommand{\om}{\omega}
\newcommand{\Pic}{\operatorname{Pic}}
\newcommand{\D}{{\mathcal D}}
\newcommand{\E}{{\mathcal E}}
\renewcommand{\O}{{\mathcal O}}
\newcommand{\C}{{\mathbb C}}
\newcommand{\Z}{{\mathbb Z}}
\newcommand{\Q}{{\mathbb Q}}
\renewcommand{\P}{{\mathbb P}}
\newcommand{\Ob}{{\mathcal O}b\,}
\newcommand{\coh}{{\operatorname{coh}}}
\newcommand{\Coh}{\operatorname{Coh}}
\newcommand{\opp}{{\operatorname{opp}}}
\newcommand{\supp}{\operatorname{supp}}
\newcommand{\Tor}{\operatorname{Tor}}
\newcommand{\Hom}{\operatorname{Hom}}
\newcommand{\Homc}{\operatorname{Hom^{\scriptscriptstyle\bullet}}}
\newcommand{\RHom}{\operatorname{RHom}}
\renewcommand{\hom}{\operatorname{\text
     {${\mathcal H}o\mskip.1\thinmuskip m$}}}
\newcommand{\Rhom}{\operatorname{\text
     {${\mathcal R}{\mathcal H}o\mskip.1\thinmuskip m$}}}
\newcommand{\h}{{\mathcal H}}
\newcommand{\modr}{{\operatorname{mod}}{-}\mskip-.1\thinmuskip}
\newcommand{\otL}{\otimes^{\mathcal L}}
\newcommand{\rarrow}{\DOTSB\longrightarrow}

\newcommand{\lrarrow}{\DOTSB\,\relbar\joinrel\relbar\joinrel\rightarrow\,}

\renewcommand{\:}{\colon}
\newcommand{\nbk}{\nobreak}

\newenvironment{thm}[1]{\smallskip{\sf #1.}\sl}{\smallskip}
\newenvironment{rem}[1]{\smallskip{\sf #1}}{\smallskip}
\newcommand{\pr}[1]{{\it #1\/}:}
\newcommand{\Qed}{\qed\smallskip}

\begin{document}
\rightline{\scriptsize\hfill
 Preprint alg-geom/9507014}
\vspace{0.75cm}
\title{All strictly exceptional collections in $\D^b_\coh(\P^m)$ \\
                     Consist of vector bundles}
\author{Leonid Positselski}
\address{Independent University of Moscow}
\email{posic@mccme.ru}
\maketitle

\section{Introduction}

 Let $\k$ be a field and $\D$ be a $\k$-linear triangulated
category; we will denote, as usually, $\Hom^i(X,Y)=\Hom(X,Y[i])$
and $\Homc(X,Y)=\bigoplus_i\Hom^i(X,Y)$.
 An object $E\in\Ob\D$ is called {\it exceptional\/} if one has
$\Hom^s(E,E)=0$ for $s\ne0$ and $\Hom^0(E,E)=\nbk\k$.
 A finite sequence $\E$ of exceptional objects $E_1,\ds,E_n$ is
called an {\it exceptional collection\/} if $\Homc(E_i,E_j)=0$ for
$i>j$.
 A collection $\E$ is called {\it full\/} if it generates $\D$ in
the sense that any object of $\D$ can be obtained from $E_i$ by
the operations of shift and cone.
 The Grothendieck group $K_0(\D)$ of a triangulated category
$\D$ generated by an exceptional collection $\E$ is the free
$\Z$-module generated by the classes of $E_1,\ds,E_n$, so
any full exceptional collection consists of $n=\rk K_0(\D)$
objects.
 Moreover, it is explained in the paper~\cite{BK} that
(under some technical restriction which is usually satisfied)
a triangulated category $\D$ generated by an exceptional collection
$\E$ is equivalent to the derived category of modules over the
differential graded algebra corresponding to $\E$.

 Let $(E_1,\,E_2)$ be an exceptional pair; the {\it left\/} and
{\it right mutated objects\/} $L_{E_1}E_2$ and $R_{E_2}E_1$ are
defined as the third vertices of exceptional triangles
 \begin{gather*}
  E_2[-1]\lrarrow L_{E_1}E_2\lrarrow\Homc(E_1,E_2)\ot E_1\lrarrow E_2 \\
  E_1\lrarrow\Homc(E_1,E_2)^*\ot E_2\lrarrow R_{E_2}E_1\lrarrow E_1[1].
 \end{gather*}
 This definition was given in the papers~\cite{Gor,Bon}; it was shown
that the {\it mutated collections\/}
 \begin{gather*}
  E_1,\ds,E_{i-2},\,L_{E_{i-1}}E_i,\,E_{i-1},\,
               E_{i+1},\ds,E_n  \\
  E_1,\ds,E_{i-2},\,E_i,\,R_{E_i}E_{i-1},\,
               E_{i+1},\ds,E_n
 \end{gather*}
remain exceptional (and full) and that the left and right
mutations are inverse to each other.
 Mutations defined this way form an action of the Artin's braid
group $B_n$ with $n$ strings on the set of all isomorphism classes
of exceptional collections of $n$ objects.

 There is a central element $\phi\in B_n$ that corresponds to the
rotation action of the circle on the space of $n$-point
configurations in $\C$.
 Its action on exceptional collections can be described as
follows.
 Let $E_{n+1}=R^{n-1}E_1$ be the object obtained by successive
left mutations of $E_1$ through $E_2$, \ds, $E_n$.
 Then it follows that the collection $E_2,\ds,E_{n+1}$ is also
exceptional.
 Proceeding in this way, we obtain the collection $E_3,\ds,E_{n+2}$,
and so on, constructing an infinite sequence of exceptional objects
$E_1$,~$E_2$,~$E_3$,~\ds with the property that any
$n$ sequential objects $E_i,\ds,E_{i+n-1}$ form an exceptional
collection.
 Using left mutations, we can continue it to the negative
indices: $E_0=L^{n-1}E_n$, $\,E_{-1}=L^{n-1}E_{n-1}$, and so on.
 This sequence is called a {\it helix}.
 The action of $\phi$ on exceptional collections shifts it
$n$ times to the left:
 $$
  \phi(E_1,\ds,E_n)=(E_{-n+1},\ds,E_0).
 $$
 The point is that this shift can be extended to an exact auto-equivalence
of the category $\D$.
 Namely, the {\it Serre functor\/} for a triangulated category
$\D$ is a covariant functor $F\:\D\rarrow\D$ for which there
is a natural isomorphism
 $$
  \Homc(U,V)=\Homc(V,FU)^*.
 $$
 It is shown in~\cite{Bon} that one has $E_{i-n}=F(E_i)[-n+1]$ for a full
helix $\E$ in $\D$.

 Now let us turn to exceptional collections in the derived
category $\D^b_\coh(X)$ of coherent sheaves on a smooth projective
algebraic variety $X$.
 In this case the Serre functor has the form
$F(U)=U\ot\om_X[\dim X]$, where $\om_X$ is the canonical line bundle.
 In the initial works of A.~Gorodentsev and A.~Rudakov~\cite{GR},
they considered exceptional collections consisting of pure sheaves,
not complexes.
 Therefore, such mutations were not defined for any exceptional
collections, but only under the conditions that some maps are
injective or surjective.
 For example, we see that the helix generated by a full exceptional
collection of sheaves will not consist of sheaves unless its period $n$
is equal to $\dim X+1$.

 Conversely, it was shown by A.~Bondal~\cite{Bon} that all mutations
of a full exceptional collection of $\dim X+1$ sheaves in
$\D^b_\coh(X)$ (that is, for a variety with $\rk K_0(X)=\dim X+1$)
consist of pure sheaves again.
 Indeed, the statement that $R_{E_2}E_1$ is a sheave follows
immediately from the isomorphism
 $R_{E_2}E_1=L_{E_3}\dsm L_{E_{n}}E_{n+1}$,
where all the objects $E_1,\ds,E_{n+1}$ are pure sheaves.
 It is also easy to see that in this case mutations preserve
the property of a full exceptional collection to consist of
locally free sheaves.
 Note that for any projective variety one has
  $\rk K_0(X)\ge\dim X+1$,
since the cycles of self-intersection of $\O(1)$ are linearly
independent over $\Q$; the equality holds for $\P^m$,
odd-dimensional quadrics, and some others.

 The principal problem of the theory of mutations of exceptional
bundles on $\P^m$ is to prove that their action on full
exceptional collections of vector bundles is transitive.
 More generally, it was conjectured in~\cite{BP} that the action
of the semidirect product group $B_n\righttimes\nbk\Z^n$ generated
by mutations and shifts on full exceptional collections in any
triangulated category $\D$ is transitive.
 The second half of this latter conjecture for smooth projective
varieties with $\rk K_0(X)=\dim X+1$ states that any full
exceptional collection in $\D^b_\coh(X)$ consists of shifts of
vector bundles.
 In this paper we prove this last statement under the following
additional restriction.
 An exceptional collection is said to be {\it strictly exceptional\/} if
one has $\Hom^s(E_i,E_j)=0$ for $s\ne0$.

\begin{thm}{Theorem}
 Let $X$ be a smooth projective variety for which
$n=\rk K_0(X)=\dim X+\nbk1$.
 Then for any strictly exceptional collection $E_1,\ds,E_n$
generating $\D^b_\coh(X)$ the objects $E_i$ are locally
free sheaves shifted on the same number $a\in\Z$ in $\D$.
 \end{thm}

 Conversely, it was shown in~\cite{Bon} that any full exceptional
collection of $\dim X+1$ sheaves on a smooth projective variety is
strictly exceptional.

 In particular, if a full exceptional collection on a variety with
$\rk K_0(X)=\dim X+1$ consists of pure sheaves, then these sheaves are
locally free.
 On the other hand, it follows that the property of a full exceptional
collection in a triangulated category of this kind to be strictly
exceptional is preserved by mutations; moreover, all strictly exceptional
collections in these categories are {\it geometric\/} in the sense
of~\cite{BP}.

At last, our methods provide an approach to the results
on recovery of algebraic varieties from the derived categories of coherent
sheaves, alternative to the one given by Bondal--Orlov~\cite{BO}.

\begin{thm}{Corollary}
 Suppose the canonical sheave of a smooth projective variety $X$
is either ample or anti-ample.
 Then the standard $t$-structure on the derived category $\D^b_{\coh}(X)$
can be recovered (uniquely up to a shift) from the triangulated category
structure.
 \end{thm}

 I am grateful to A.~Bondal who introduced me into the subject of
triangulated categories and exceptional collections and to
A.~Polishchuk and A.~N.~Rudakov for very helpful discussions.
 I am pleased to thank Harvard University for its hospitality during
preparation of this paper.

\smallskip
\section{Reduction to a Local Problem}

The next result is due to A.~Bondal and A.~Polishchuk~\cite{BP}.

\begin{thm}{Proposition}
 Suppose a helix $\{E_i,\,i\in\Z\}$ in a triangulated category $\D$
is generated by a strictly exceptional collection $E_1,\ds,E_n$.
 Then one has $\Hom^s(E_i,E_j)=0$ for $s>0$ and $i\le j\in\Z$,
as well as for $s<n-1$ and $i\ge j\in\Z$.
 \end{thm}

\pr{Proof}
 First note that the Serre duality isomorphisms
 $$
  \Hom^s(E_i,E_j)=\Hom^{n-1-s}(E_{j+n},E_i)^*
 $$
mean that two statements are equivalent to each other;
let us prove the first one.
 The simplest way is to identify $\D$ with the
derived category of modules over the homomorphism algebra
$A=\bigoplus_{k,l=1}^nA_{kl}$, $\,A_{kl}=\Hom(E_k,E_l)$ of our strictly
exceptional collection, so that the objects $E_l$ correspond to the
projective $A$-modules $P_l=\bigoplus_kA_{kl}$ for $1\le l\le n$.
 Since the Serre functor provides $n$-periodicity isomorphisms
$\Hom^s(E_i,E_j)=\Hom^s(E_{i+n},E_{j+n})$, we can assume that
$1\le i\le n$.
 Let $j=k+Nn$ for some $1\le k\le n$; then we have
$E_j=F^{-N}E_l\mskip.7\thinmuskip[N(n-1)]$.
 The Serre functor on $\D^b(\modr A)$ has the form
$F(M)=\Hom_\k(\RHom_A(M,A),\k)$ and
$F^{-1}(M)=\RHom_A(\Hom_\k(M,\k),A)$;
since the homological dimension of $A$ is not greater than $n-1$,
we obtain $E_j\in\D^{\le0}(\modr {A})$ for $j\ge1$.
 Since $E_i$ are projective for $1\le i\le n$, the assertion follows.
 A direct, but more complicated calculation from~\cite{BP} allows to avoid
the additional condition on $\D$.
 \Qed

\pr{Proof of Theorem}
 First let us show that $X$ is a Fano variety.
 We give a simple strengthening of the argument from~\cite{BP}.
 Since $\rk\Pic(X)=1$, there are only three types of invertible
sheaves: ample ones, antiample ones, and sheaves of finite order.
 We have to prove that $\om^{-1}$ is ample; it is enough to show that
$H^0(\om^N)=0$ for all $N>0$.
 Let us denote by $\h^s(U)$ the cohomology sheaves of a complex $U$.
 Since $E_1,\ds,E_n$ generate $\D$, it is clear that there exists
$i$ and $s$ such that $\supp \h^s(E_i)=X$.
 Let we have a nonzero section $f\in H^0(\om^N)$; it induces
a morphism $E_i\rarrow E_i\ot\om^N$ which is nonzero since
its restriction to $\h^s$ is.
 But we have $E_i\ot\om^N=E_{i-Nn}$ which provides a
contradiction with Proposition.

\begin{rem}{Remark 1:}
 More generally, one can see that the canonical sheave $\om$ cannot
be of finite order for a variety $X$ admitting a full exceptional
collection.
 Indeed, the action of invertible sheaves on $K_0(X)$ is unipotent
with respect to the filtration by the dimensions of supports, thus
in the case in question the action of $\om$ on $K_0(X)$ must be
trivial.
 But this action (skew-)symmetrizes the canonical bilinear form
$\chi([U],[V])=\sum(-1)^s\dim\Hom^s(U,V)$ on $K_0(X)$.
 In the basis of $K_0$ corresponding to an exceptional collection,
the matrix of this form is upper-triangular with units on the
diagonal, so it cannot be skew-symmetric and if it is symmetric
then it is positive.
 The latter is impossible since one has $\chi([\O_x],[\O_x])=0$
for the structure sheave $\O_x$ of a point $x\in X$.
\end{rem}

 We will essentially use the tensor structure on $\D^b_\coh(X)$.
 Namely, let
  $$
   \Rhom\:\D^\opp\times\D\rarrow\D
  $$
be the derived functor of local homomorphisms of coherent sheaves;
it can be calculated using finite locally free resolvents.
 We have $\Hom^s(U,V)=H^s(\Rhom(U,V))$, where $H$ denotes the
global sheave's cohomology.
 Let $i$, $j\in\Z$ be fixed and $N$ be large enough; one has
  $$
   \Hom^s(E_i,E_{j+Nn})=H^s(\Rhom(E_i,E_j\ot\om^{-N}))
   =H^s(C_{ij}\ot\om^{-N}),
  $$
where we denote $C_{ij}=\Rhom(E_i,E_j)$.
 Since $\om^{-1}$ is ample, for large $N$ we have
$H^{>0}\h^s(C_{ij}\ot\om^{-N})=0$, hence by the spectral
sequence
 $$
  H^s(C_{ij}\ot\om^{-N})=H^0\h^s(C_{ij}\ot\nbk\om^{-N}).
 $$
 Using the property of ample sheaves again, we see that
$\Hom^s(E_i,E_{j+Nn})$ is nonzero iff $\h^s(C_{ij})$ is.
 Let $\D^{\le0}$ and $\D^{\ge0}$ denote the subcategories
of $\D^b_\coh(X)$ defined in the standard way.
 Comparing the last result with Proposition, we finally obtain
$C_{ij}\in\D^{\le0}$.

\begin{rem}{Remark 2:}
 Now we can show easily that our exceptional collection is
{\it geometric\/}~\cite{BP}.
 Indeed, using the duality $\Rhom(V,U)=\Rhom(\Rhom(U,V),\,\O)$,
one obtains $C_{ij}=\Rhom(C_{ji},\O)$ and since
$\Rhom(\D^{\le0},\O)\subset\D^{\ge0}$, it follows that $C_{ij}$
are pure sheaves.
 Therefore $\Hom^{<0}(E_i,E_j)=H^{<0}(C_{ij})=0$ for any $i$ and
$j$.
\end{rem}

 The following local statement allows to finish the proof.

\begin{thm}{Main Lemma}
 Let $E\in\D^b_\coh(X)$ be a coherent complex on a smooth
algebraic variety $X$ such that $\Rhom(E,E)\in\D^{\le0}$.
 Then $E$ is a (possibly shifted) locally free sheave.
  \end{thm}

 It only remains to show that all of $E_i$ are placed in the same degree in
$\D$, which is true since they are locally free and $\Rhom(E_i,E_j)$ is
placed in degree~0.  \Qed

\pr{Proof of Corollary}
 The functor of twisting on $\om$ on the derived category
can be recovered in terms of the Serre functor.
 Let $\om$ be anti-ample.
 According to Main Lemma, an object $E\in\D^b_\coh(X)$ is a
shifted vector bundle iff $\Hom^s(E,E\ot\om^{-N})=0$ for
$s\ne0$ and $N$ large enough.
 For a nonzero vector bundle $E$ and $U\in\D$ one has
$U\in\D^{\ge0}$ iff $\Hom^{<0}(E_i,U\ot\om^{-N})=0$ for
large $N$, and the same for $\D^{\le0}$.
\qed

\smallskip
\section{The Proof of Main Lemma}

\begin{thm}{Lemma 1}
 If $E\in\Coh(X)$ and $\Rhom(E,\O)$ are pure sheaves placed in
degree~0, then $E$ is locally free.
 \end{thm}

\pr{Proof}
 Let $0\rarrow P_k\rarrow P_{k-1}\rarrow\ds\rarrow P_0\rarrow0$ be a
locally free resolvent of $E$.
 Since $\hom^k(E,\O)=0$, we see that the morphism
$\hom(P_{k-1},\O)\rarrow\hom(P_k,\O)$ is surjective.
 Thus, the inclusion $P_k\rarrow P_{k-1}$ is locally split and
the quotient sheave $P_{k-1}/P_k$ is locally free, which allows to
change our resolvent to a shorter one. \Qed

 Let $U\otL V$ denote the derived functor of tensor product
over $\O_X$ on $\D^b_\coh(X)$; then one has
$\Rhom(U,V)=\Rhom(U,\O)\otL V$.

\begin{thm}{Lemma 2}
 Let $E$, $F\in\D^b_\coh(X)$; suppose $E\otL F\in\D^{\le0}$.
 Then for any $i+j\ge0$ one has
$\supp\h^i(E)\cap\supp\h^j(F)=\empty$.
 \end{thm}

\pr{Proof}
 Proceed by decreasing induction on $i+j$.
 Consider the K\"unneth spectral sequence
  $$
   E_2^{pq}=\bigoplus_{i+j=q}\Tor_{-p}(\h^iE,\h^jF)
   \implies \h^{p+q}(E\otL F).
  $$
 If the intersection of supports is nonzero, then it is easy to
see that $\h^iE\ot\h^jF\ne0$, thus $E_2^{0,q}\ne0$.
 This term can be only killed by some
$E^{-r,q+r-1}$, where $r\ge2$; but it follows from the induction
hypothesis that $E_2^{p,\ge q+1}=0$. \Qed

\pr{Proof of Main Lemma}
 Let $F=\Rhom(E,\O)$; then one has $\Rhom(E,E)=E\otL F$.
 Using a shift, we can assume that $E\in\D^{\le0}$ and $\h^0(E)\ne0$;
then $F\in\D^{\ge0}$ and $\h^0F=\hom(\h^0E,\O)$.
 By Lemma~2, we have $\supp\h^0E\cap\supp\h^{>0}F=\empty$.
 Clearly, one can assume that $X$ is irreducible.
 First let us show that $\supp\h^0(E)=X$.
 Indeed, in the other case it is clear that $\h^0F=0$ and the
restriction of $F$ on $X\setminus\supp\h^{>0}F$ is acyclic while
the restriction of $E$ is not, which contradicts the local nature
of $\Rhom$.
 Thus we have $\supp\h^0(E)=X$, which implies $\h^{>0}F=0$ and
$F\in\Coh(X)$.
 It follows that $E=\Rhom(F,\O)\in\D^{\ge0}$ and $E\in\Coh(X)$.
 By Lemma~1, $E$ is locally free. \Qed

\end{document}